\documentclass[11pt]{amsart}
\usepackage{amsmath,amssymb,amsthm,amscd,latexsym}
\usepackage[all]{xy}
\usepackage{graphicx}
\usepackage{hyperref}

\input xypic
\xyoption{all}

\newcommand{\rar}{\rightarrow}
\newcommand{\lar}{\longrightarrow}
\newcommand{\llar}{-\kern-5pt-\kern-5pt\longrightarrow}

\newtheorem{Theorem}{Theorem}[section]
\newtheorem{Lemma}[Theorem]{Lemma}
\newtheorem{Corollary}[Theorem]{Corollary}
\newtheorem{Proposition}[Theorem]{Proposition}
\newtheorem{Remark}[Theorem]{Remark}
\newtheorem{Example}[Theorem]{Example}
\newtheorem{Conjecture}[Theorem]{Conjecture}

\def\sqr#1#2{{\vcenter{\hrule height.#2pt
        \hbox{\vrule width.#2pt height#1pt \kern#1pt
            \vrule width.#2pt}
        \hrule height.#2pt}}}
\def\phi{\varphi}
\def\demo{\noindent{\bf Proof. }}
\def\square{\mathchoice\sqr64\sqr64\sqr{4}3\sqr{3}3}
\def\qed{\hspace*{\fill} $\square$}

\def\fm{{\mathfrak m}}

\def\Ree#1{{\mathcal R}(#1)}

\def\restr{{\kern-1pt\restriction\kern-1pt}}

\def\pp{{\mathbb P}}

\begin{document}
\setlength{\hoffset}{-2cm}
\title[Binary general forms]{The depth of the Rees algebra of three general binary forms}

\author{Ricardo Burity}


\author{Aron Simis}

\subjclass[2010]{13A02, 13A30, 13C15, 13D02, 13D40.} \keywords{Rees algebra,  almost Cohen--Macaulay, Hilbert function, Ratliff--Rush filtration, Huckaba--Marley test. \\
\indent The second author was partially supported by a CNPq grant (302298/2014-2) and various visiting fellowships.\\
\indent Burity address: Departamento de Matem\'atica, Universidade Federal da Paraiba, 58051-900 J. Pessoa, PB, Brazil, email: ricardo@mat.ufpb.br\\
\indent Simis address: Departamento de Matem\'atica, Universidade Federal de Pernambuco, 50740-540 Recife, PE, Brazil, email: aron@dmat.ufpe.br}


\begin{abstract}

One proves that the Rees algebra of an ideal generated by three general binary forms of same degree $\geq 5$ has depth one. The proof hinges on the behavior of the Ratliff--Rush filtration for low powers of the ideal and on establishing that certain large matrices whose entries are quadratic forms have maximal rank. One also conjectures a shorter result that implies the main theorem of the paper.

\end{abstract}

\maketitle


\section{Introduction}

Set $R:=k[x,y]$ denote a standard graded polynomial ring in two variables over an infinite field $k$. 
Let $I\subset R$ stand for a codimension $2$ ideal
generated by $3$ forms of the same degree $d\geq 2$.
The goal of this paper is to prove the following result:

\medskip

{\bf Theorem.}
Let $I=\subset R:=k[x,y]$ denote a codimension $2$ ideal generated by three general forms $f_1,f_2,f_3$ of degree $d\geq 5$.
Then the Rees algebra of $I$ has depth $1$.

\medskip

The Rees algebra of $I$ is the standard graded $R$-algebra $\Ree I: =\bigoplus_{\ell\geq 0} I^{\ell}$. Its graded maximal ideal is $\mathcal{M}:=(x,y, \bigoplus_{\ell\geq 1} I^{\ell})$.
The depth of $\Ree I$ is computed with respect to $\mathcal{M}$.
The theorem says that this depth is smallest possible.

We observe that the Rees algebra of an $(x,y)$-primary ideal of $k[x,y]$ generated by  three arbitrary forms of degree $d\leq 4$ is almost Cohen--Macaulay: when $I$ admits a linear syzygy, the result is part of \cite[Theorem 4.4]{syl3} (see also \cite[Theorem 4.4]{KPU} and \cite[Proposition 2.3]{ST}).
The balanced case $d=4$ is proved in \cite[Proposition 4.3]{syl3}.

The above theorem is in addition in sharp contrast to recent akin statements regarding the depth of $\Ree I$ when $I$ is an almost complete intersection (see \cite{syl1,syl3,RoSw,ST}), where it has been proved that the depth is at least $2$ -- i.e.,  $\Ree I$ is almost Cohen--Macaulay. 

Yet another difference is that the present methods hardly touch directly the structure of the presentation ideal of $\Ree I$ as an $R$-algebra. In fact, the entire matter is pretty much decided at the level of the second and third powers of the ideal $I$ through the use of two apparently disconnected tools: the Ratliff-Rush filtration and the Huckaba--Marley test.

In order to use the first of these tools we were lead to ``solve'' large linear systems over the ground field, and hence to decide whether large matrices have maximal rank. The surest way to go about that was to go all the way to some generic situation. The price to pay is the reader's indulgence in following some large but elementary discussion on matrices whose entries are quadratic forms in many variables with unit coefficients.

A couple of words about the terminology ``general forms''.
One usually says that a set of forms in a polynomial ring over an infinite field is {\em general} in the sense that the total collection of the coefficients of the forms is general in the parameter space of the coefficients.
In a concrete situation, this is understood in the sense that one avoids a contextualized closed set in the parameter space.
Such a property is often hard to work with due to its instability under ordinary algebraic operations. Even the common perception that over three or more variables a general form is irreducible over $k$ becomes elusive in the present case of two variables since in this case any form factors into linear forms over an algebraic closure of $k$. 
Fortunately, for computation purposes taking random coefficients throughout will do. 

\medskip

{\sc Acknowledgements:} The authors thank Stefan Toh\v aneanu for early discussions in many passages and basic calculations, specially in the first two sections.

\section{Preliminaries on three binary general forms}

Set $R:=k[x,y]$, a standard graded polynomial ring in two variables over an infinite field $k$. 
Let $I\subset R$ denote a codimension $2$ ideal
generated by $3$ forms of the same degree $d\geq 2$.
We consider a minimal free resolution of $I$
\begin{equation}\label{HB2}
0\rar R(-(d+r))+R(-(d+s))\stackrel{\phi}{\lar} R(-d)^3\lar I\rar 0,
\end{equation}
where $1\leq s\leq r$ denote the standard degrees of the columns of a matrix of $\phi$.
We observe that $I$ is the ideal generated by the $2\times 2$ minors of $\phi$, and hence
one often says that $\phi$ is the {\em Hilbert--Burch matrix associated} to $I$ and $s,r$ are its standard degrees.
Note that the $r,s$ are numerical invariants of $I$ adding up to $d$.

\begin{Lemma}\label{semi-balanced}
Let $I\subset R$ denote a codimension $2$ ideal
generated by $3$ general forms of the same degree $d\geq 2$.
Then the syzygy module of $I$ is generated in degrees
 $d+\left \lfloor{d/2}\right \rfloor$ and $d+\left \lceil{d/2}\right \rceil$.
\end{Lemma}
\demo 
We induct on $d$.
Since $d$ may be even or odd and our inductive process will pass from $d-2$ to $d$, we need to start the induction in the cases where $d=2$ and $d=3$, respectively.
Obviously, there is nothing to show in these two cases as the statement holds by default for any set of forms in these degrees.

We now proceed to the inductive step, assuming $d\geq 4$ for the even case and $d\geq 5$ for the odd case.
Up to general $k$-linear combinations we can assume that the given forms have the shape
\begin{eqnarray}\label{triangular}
f_1&=&x^d+y^3g_1\nonumber\\
f_2&=&x^{d-1}y+y^3g_2 \\
f_3&=&x^{d-2}y^2+y^3g_3.\nonumber
\end{eqnarray}
Note that the set  $\{g_1,g_2,g_3\}$ of forms preserve almost all coefficients of the original $f$'s, hence is general. Their common degree is $d-3$.

Suppose $Pf_1+Qf_2+Rf_3=0$ is a syzygy. Since $f_2,f_3$ are divisible by $y$, one has $\boxed{P=yP'}$, for some $P'\in k[x,y]$. Plugging this in the syzygy equation, and canceling  $y$ yields
$$P'(x^d+y^3g_1)+Q(x^{d-1}+y^2g_2)+R(x^{d-2}y+y^2g_3)=0.$$
Then $y$ divides $x^dP'+Qx^{d-1}=x^{d-1}(Q+xP')$, giving $\boxed{Q=yQ'-xP'}$ for some $Q'\in k[x,y]$. Plugging this back into the above equation, after canceling $P'x^d$, and after simplifying by $y$, one has
$$P'y^2g_1+Q'x^{d-1}+y^2Q'g_2-xyP'g_2+Rx^{d-2}+yRg_3=0.$$
Thus, $y$ divides $Q'x^{d-1}+Rx^{d-2}=x^{d-2}(R+xQ')$, giving $\boxed{R=yR'-xQ'}$, for some $R'\in k[x,y]$.

Again, plugging this back, canceling $Q'x^{d-1}$, and simplifying by $y$ yields
$$P'(yg_1-xg_2)+Q'(yg_2-xg_3)+R'(x^{d-2}+yg_3)=0.$$

Note that the coefficients of the three forms $\{yg_1-xg_2,yg_2-xg_3,x^{d-2}+yg_3\}$  of degree $d-2$ are almost all of the original general coefficients or sums of these, hence this set is again a general set of forms of degrees $d-2$.
By the inductive hypothesis on $d$ they generate an ideal whose  syzygy module is generated in standard degrees $\left \lfloor{(d-2)/2}\right \rfloor$ and $\left \lceil{(d-2)/2}\right \rceil$.

Using the formulas in the boxes, we get two independent syzygies on $(f_1,f_2,f_3)$ of standard degrees $\left \lfloor{d/2}\right \rfloor$ and $\left \lceil{d/2}\right \rceil$. \qed

\medskip

Let us read  the Hilbert function $H(I,t)$ of $I$ out of Lemma~\ref{semi-balanced}:

\begin{Lemma}\label{Hilbert_function}
Let $I\subset R$ denote a codimension $2$ ideal
generated by $3$ general forms of the same degree $d\geq 2$.
Then
$$H(I,t)=3H(R,t-d)-H(R,t-d-\left \lfloor{d/2}\right \rfloor)-H(R,t-d-\left \lceil{d/2}\right \rceil).$$
Consequently, setting $s:=\left\lfloor{d/2}\right \rfloor$, one has:
\begin{enumerate}
\item[{\rm (Odd)}] If $d$ is odd then  $\dim (R/I)_{3s}=1$ and $\dim (R/I)_{3s+1}=0$.
\item[{\rm (Even)}] If $d$ is even then  $\dim (R/I)_{3s-2}=2$ and $\dim (R/I)_{3s-1}=0$.
\end{enumerate}

\noindent Moreover, up to $k$-linear transformations and change of variables, we can assume the following:

$\bullet$ If $d=2s+1$, then $y^{3s}$ spans $(R/I)_{3s}$; 

$\bullet$  If $d=2s$, then $\{y^{3s-2}, xy^{3s-3}\}$ span $(R/I)_{3s-2}$.
\end{Lemma}
\demo
The Hilbert function easily reads off  the minimal graded resolution
$$0\rar R(-d-\left\lfloor{d/2}\right \rfloor)\oplus R(-d-\left\lceil{d/2}\right \rfloor)\lar R(-d)^3 \lar I\rar 0$$
afforded by Lemma~\ref{semi-balanced}.

Therefore, the full Hilbert function of $R/I$ is
$$\begin{array}{ccccccc}
{\rm dimension}\quad\quad\quad \ldots 2s+1 & 2s-1 & 2s-3 & \ldots &  3  &  1 &  0 \\
{\rm degree}\quad\quad\quad\quad\;\; \ldots    2s & 2s+1 & 2s+2 & \ldots & 3s-1 & 3s & 3s+1
\end{array},
$$
if $d=2s+1$, and

$$\begin{array}{ccccccc}
{\rm dimension}\quad\quad\quad \ldots 2s & 2s-2 & 2s-4 & \ldots &  4  &  2 &  0 \\
\kern0.7cm{\rm degree}\quad\quad\quad  \ldots 2s-1 & 2s & 2s+1 & \ldots & 3s-3 & 3s-2 & 3s-1
\end{array},
$$
if $d=2s$.

\smallskip

In order to prove the supplementary assertion, we focus on the odd case $d=2s+1$,  the even case being similar.
Consider the $k$-vector space $I_{3s}=(x,y)_{s-1}I_{2s+1}$ and recall that it has diimension $3s$.
By a similar procedure as (\ref{triangular}) in the proof of Lemma~\ref{semi-balanced}, we can apply a suitable action of $GL(3s,k)$ to obtaining a $k$-basis of $I_{3s}$ in a ``triangular'' form 

\begin{eqnarray}
F_1&=&\alpha_{3s,0}x^{3s}+ \alpha_{3s-1,1}x^{3s-1}y+\cdots\nonumber\\
F_2&=&\quad\quad\quad\quad \beta_{3s-1,1}x^{3s-1}y+ \beta_{3s-2,2}x^{3s-2}y^2+\cdots\nonumber\\
F_3&=&\quad\quad\quad\quad\quad\quad\quad\quad\quad\quad \gamma_{3s-2,2}x^{3s-3}y^2+\beta_{3s-3,3}x^{3s-3}y^3+\cdots\nonumber\\
&\vdots& \quad\quad\quad\quad \quad\quad\quad\quad\quad\quad\quad\quad\quad\quad\vdots\nonumber\\
F_{3s-1}\kern-6pt&=&\quad\quad\quad\quad\quad\quad\quad\quad\quad\quad\quad\quad\quad\mu_{2,3s-2}x^2y^{3s-2}+\mu_{1,3s-1}xy^{3s-1}+\mu_{3s}y^{3s}\nonumber\\
F_{3s}&=& \quad\quad\quad\quad\quad\quad\quad\quad\quad\quad\quad\quad\quad\quad\quad\quad\quad\quad\quad\quad\nu_{1,3s-1}xy^{3s-1} +\nu_{3s} y^{3s}\nonumber
\end{eqnarray}

Applying the change of variables 
$$x\mapsto \frac{1}{\nu_{1,3s-1}}x-\frac{\nu_{3s} }{\nu_{1,3s-1}},\;\; y\mapsto y
$$
transforms $F_{3s}$ in the above basis into $xy^{3s-1}$.
Let $\tilde{I}\subset R$ denote the ideal obtained by applying this change of variables to the forms generating $I$.
Then, since $I_{3s}=(x,y)_{s-1}I_{2s+1}$ and $(x,y)$ is invariant under this change, it follows that $\tilde{I}_{3s}=(x,y)_{s-1}\tilde{I}_{2s+1}$ contains the monomial $xy^{3s-1}$. 
On the other hand, the change of variables did not affect the term $\mu_{3s}y^{3s}$ of the basis element $F_{3s-1}$.
If $y^{3s}\in \tilde{I}$ then we can clean bottom up all monomials, thus implying that $I_{3s}$ is spanned by all monomials of this degree - an absurd. 
\qed

\section{A conjectured predecessor}

For the conjectural statement of this section we focus on the case of odd degree $d=2s+1$. The case of even degree may admit a similar treatment.

We will need to know that in the Hilbert function above any monomial of degree $3s$ spans $(R/I)_{3s}$ -- the proof is the same as in Lemma~\ref{Hilbert_function} for the monomail $y^{3s}$.

We establish the technical preliminaries.
Given a form 
\begin{eqnarray*}\nonumber
f&=&\alpha_{d,0}x^d+\alpha_{d-1,1}x^{d-1}y+\cdots + \alpha_{2,d-1}x^2y^{d-2}+\alpha_{1,d-1}xy^{d-1}+\alpha_{0,d}y^d \\ \nonumber
&=& (\alpha_{d,0}x^{d-2}+\alpha_{d-1,1}x^{d-3}y+\cdots + \alpha_{2,d-1}y^{d-2}) \,x^2+\alpha_{1,d-1}xy^{d-1}+\alpha_{0,d}y^d
\end{eqnarray*}
of degree $d$ we associate to it the unique form $\delta(f)$ of degree $d-2$ which is the coefficient of $x^2$ in the above expression.
It is clear that this association preserves many properties of the original form, including that of being general provided $d >> 0$.
Now, given three general forms $f_1,f_2,f_3$ of degree $d$, in the sense that the $(d+1)\times 3$ total matrix of the coefficients has general entries, we let
$I=(f_1,f_2,f_3)\subset R$ and  $\delta(I):=(\delta(f_1),\delta(f_2), \delta(f_3))\subset R$.

By construction, $I\subset \widetilde{I}:=(x^2\delta(I), xy^{d-1}, y^d)$; more precisely, this is an inclusion induced by an inclusion of the corresponding linear systems spanned in the common degree $d$.

Set $\fm:=(R_+)=(x,y)$.

\begin{Conjecture} Let $d=2s+1\geq 7$. Then $\fm^{s-1}(\widetilde{I})^2\subset I^2$.
\end{Conjecture} 

\begin{Remark}\rm
It may be the case that actually $I^2:(\widetilde{I})^2=\fm^{t}$, for some $t\leq s-1$ but for the purpose to follow the above inclusion is sufficient.
\end{Remark} 

\begin{Corollary} {\rm (To the conjecture)}
	Let $d=2s+1\geq 5$. Then the annihilator of the $R/I$-module $I/I^2$ contains a monomial spanning $(R/I)_{3s}$, and hence $I^2:I\neq I$.
\end{Corollary}
\demo We induct on $s\geq 2$.
The initial case $s=2$ is to be verified directly by computation (e.g., by the method in the last section). 
 
Let then $s\geq 3$.
By the inductive hypothesis, we may assume that $y^{3(s-1)}\in \delta(I)^2:\delta(I)$.
More precisely, by a degree count one has $y^{3(s-1)}\delta(I)\subset \fm^{s-2}\delta(I)^2$.

We have noted above that regarding the Hilbert function of $R/I$, $x^2y^{3s-2}$ spans $(R/I)_{3s}$.
For this monomial one has:

$$x^2y^{3s-2}I\subset x^2y^{3s-2}\widetilde{I}=x^2y^{3s-2}(x^2\delta(I),xy^{2s},y^{2s+1}).$$

Treat each piece separately, using the conjectured result:

(i) $x^2y^{3s-2}x^2\delta(I)= x^4y(y^{3s-3}\delta(I))\subset x^4y(\fm^{s-2}\delta(I)^2)
=ym^{s-2}(x^2J)^2\subset ym^{s-2}(\widetilde{I})^2\subset \fm^{s-1}(\widetilde{I})^2\subset I^2$.

\smallskip

(ii) $x^2y^{3s-2}\, xy^{2s}= xy^{s-2}(x^2y^{4s})\subset \fm^{s-1}(\widetilde{I})^2\subset I^2$.

\smallskip

(iii) $x^2y^{3s-2} \, y^{2s+1}= x^2y^{s-3}(y^{2(2s+1)})\subset \fm^{s-1}(\widetilde{I})^2\subset I^2$.
\qed


As we will argue in the subsequent sections this conjectured statement suffices for the purpose the main result of the paper.

\section{A lemma on large matrices}

This part is about the rank of a particular type of matrix whose entries are either zeros or certain quadrics in a given set of variables over an infinite field (of characteristic $\neq 2$). As usual while dealing with such matrices, settling an argument will be shorter than transcribing the details of the matrices themselves.

\begin{Lemma}\label{main_lemma} Let $k$ denote an infinite field and let $d$ be a positive integer. 
	\begin{enumerate}
		\item[{\rm (a)}] If $d$ is odd, assume that $d\geq 5$ and write $d=2s+1$.
	Let 
	$$S:=k[T_l^{(t)}\,|\, 1\leq t\leq 3\,; 0\leq l\leq d=2s+1]$$ 
stand for a polynomial ring in $3(d+1)=6(s+1)$ variables. 
For each pair of indices $1\leq t_1\leq t_2\leq 3$, introduce the $2$-forms

$$
Q_r^{t_1,t_2} := \left\{
\begin{array}{rcl}
\displaystyle\sum_{l_1=0}^r T_{l_1}^{(t_1)}T_{r-l_1}^{(t_2)}, & \mbox{if} & 0\leq r\leq d=2s+1\\
\displaystyle\sum_{l_1=r-(2s+1)}^{2s+1} T_{l_1}^{(t_1)}T_{r-l_1}^{(t_2)}, & \mbox{if} & 2s+2\leq r\leq 2d=4s+2
\end{array}
\right.
$$

and consider the following $(5s+2)\times s$ matrix
{\small
$$B_{t_1,t_2}=\left[\begin{array}{cccccc}
Q_0^{t_1,t_2}& 0& 0& \cdots& 0 & 0 \\ [2pt]
Q_1^{t_1,t_2}& Q_0^{t_1,t_2}& 0& \cdots& 0 & 0  \\ [2pt]
Q_2^{t_1,t_2}&Q_1^{t_1,t_2}&Q_0^{t_1,t_2}&\cdots & 0& 0 \\ [2pt]
\vdots&&&&\vdots&\vdots \\
Q_{s-2}^{t_1,t_2}&Q_{s-3}^{t_1,t_2}&Q_{s-4}^{t_1,t_2}&\cdots & Q_0^{t_1,t_2}& 0 \\ [4pt]
Q_{s-1}^{t_1,t_2}&Q_{s-2}^{t_1,t_2} &Q_{s-3}^{t_1,t_2} & \cdots& Q_1^{t_1,t_2}&Q_0^{t_1,t_2} \\ [4pt]
Q_{s}^{t_1,t_2}&Q_{s-1}^{t_1,t_2} &Q_{s-2}^{t_1,t_2} & \cdots&Q_2^{t_1,t_2}& Q_1^{t_1,t_2}  \\
\vdots&&&&\vdots&\vdots \\
Q_{2s+1}^{t_1,t_2}&Q_{2s}^{t_1,t_2} &Q_{2s-1}^{t_1,t_2} & \cdots& Q_{s+3}^{t_1,t_2}&Q_{s+2}^{t_1,t_2} \\ [4pt]
Q_{2s+2}^{t_1,t_2}&Q_{2s+1}^{t_1,t_2} &Q_{2s}^{t_1,t_2} & \cdots&Q_{s+4}^{t_1,t_2} &Q_{s+3}^{t_1,t_2}  \\
\vdots&&&&\vdots&\vdots \\ 
Q_{4s+1}^{t_1,t_2}&Q_{4s}^{t_1,t_2} &Q_{4s-1}^{t_1,t_2} & \cdots& Q_{3s+3}^{t_1,t_2}& Q_{3s+2}^{t_1,t_2} \\ [4pt]
Q_{4s+2}^{t_1,t_2}&Q_{4s+1}^{t_1,t_2} &Q_{4s}^{t_1,t_2} & \cdots&Q_{3s+4}^{t_1,t_2} &Q_{3s+3}^{t_1,t_2} \\ [4pt]
0&Q_{4s+2}^{t_1,t_2} &Q_{4s+1}^{t_1,t_2} & \cdots& Q_{3s+ 5}^{t_1,t_2} & Q_{3s+4}^{t_1,t_2} \\ [4pt]
0 & 0 &Q_{4s+2}^{t_1,t_2} & \cdots& Q_{3s+6}^{t_1,t_2}&Q_{3s+5}^{t_1,t_2}\\ [4pt]
0 & 0 & 0 & \cdots& Q_{3s+7}^{t_1,t_2}&Q_{3s+6}^{t_1,t_2} \\
\vdots&&&&\vdots&\vdots \\
0 & 0 & 0&\cdots& Q_{4s+2}^{t_1,t_2}& Q_{4s+1}^{t_1,t_2} \\ [4pt]
0 & 0 & 0&\cdots& 0& Q_{4s+2}^{t_1,t_2}  \\
\end{array}\right]
$$}
where every $2$-form above is orderly and diagonally displayed $s$ times.  
Then the concatenation
$$A:=[B_{1,1}~B_{1,2}~B_{1,3}~B_{2,2}~B_{2,3}~B_{3,3}]$$
is a $(5s+2)\times 6s$ matrix of maximal rank.

\item[{\rm (b)}] If $d$ is even and $d\geq 10$, write $d=2s$.
Let 
$$S:=k[T_l^{(t)}\,|\, 1\leq t\leq 3\,; 0\leq l\leq d=2s]$$ 
stand for a polynomial ring in $3(d+1)=3(2s+1)$ variables. 
For each pair of indices $1\leq t_1\leq t_2\leq 3$, introduce the $2$-forms

$$
Q_r^{t_1,t_2} := \left\{
\begin{array}{rcl}
\displaystyle\sum_{l_1=0}^r T_{l_1}^{(t_1)}T_{r-l_1}^{(t_2)}, & \mbox{if} & 0\leq r\leq d=2s\\
\displaystyle\sum_{l_1=r-2s}^{2s} T_{l_1}^{(t_1)}T_{r-l_1}^{(t_2)}, & \mbox{if} & 2s+1\leq r\leq 2d=4s
\end{array}
\right.
$$

and consider the following $(5s-1)\times (s-1)$ matrix
{\small
	$$B_{t_1,t_2}=\left[\begin{array}{cccccc}
	Q_0^{t_1,t_2}& 0& 0& \cdots& 0 & 0 \\ [2pt]
	Q_1^{t_1,t_2}& Q_0^{t_1,t_2}& 0& \cdots& 0 & 0  \\ [2pt]
	Q_2^{t_1,t_2}&Q_1^{t_1,t_2}&Q_0^{t_1,t_2}&\cdots & 0& 0 \\ [2pt]
	\vdots&&&&\vdots&\vdots \\
	Q_{s-3}^{t_1,t_2}&Q_{s-4}^{t_1,t_2}&Q_{s-5}^{t_1,t_2}&\cdots & Q_0^{t_1,t_2}& 0 \\ [4pt]
	Q_{s-2}^{t_1,t_2}&Q_{s-3}^{t_1,t_2} &Q_{s-4}^{t_1,t_2} & \cdots& Q_1^{t_1,t_2}&Q_0^{t_1,t_2} \\ [4pt]
	Q_{s-1}^{t_1,t_2}&Q_{s-2}^{t_1,t_2} &Q_{s-3}^{t_1,t_2} & \cdots&Q_2^{t_1,t_2}& Q_1^{t_1,t_2}  \\
	\vdots&&&&\vdots&\vdots \\
	Q_{2s}^{t_1,t_2}&Q_{2s-1}^{t_1,t_2} &Q_{2s-2}^{t_1,t_2} & \cdots& Q_{s+3}^{t_1,t_2}&Q_{s+2}^{t_1,t_2} \\ [4pt]
	Q_{2s+1}^{t_1,t_2}&Q_{2s}^{t_1,t_2} &Q_{2s-1}^{t_1,t_2} & \cdots&Q_{s+4}^{t_1,t_2} &Q_{s+3}^{t_1,t_2}  \\
	\vdots&&&&\vdots&\vdots \\ 
	Q_{4s-1}^{t_1,t_2}&Q_{4s-2}^{t_1,t_2} &Q_{4s-3}^{t_1,t_2} & \cdots& Q_{3s+2}^{t_1,t_2}& Q_{3s+1}^{t_1,t_2} \\ [4pt]
	Q_{4s}^{t_1,t_2}&Q_{4s-1}^{t_1,t_2} &Q_{4s-2}^{t_1,t_2} & \cdots&Q_{3s+3}^{t_1,t_2} &Q_{3s+2}^{t_1,t_2} \\ [4pt]
	0&Q_{4s}^{t_1,t_2} &Q_{4s-1}^{t_1,t_2} & \cdots& Q_{3s+ 4}^{t_1,t_2} & Q_{3s+3}^{t_1,t_2} \\ [4pt]
	0 & 0 &Q_{4s}^{t_1,t_2} & \cdots& Q_{3s+5}^{t_1,t_2}&Q_{3s+4}^{t_1,t_2}\\ [4pt]
	0 & 0 & 0 & \cdots& Q_{3s+6}^{t_1,t_2}&Q_{3s+5}^{t_1,t_2} \\
	\vdots&&&&\vdots&\vdots \\
	0 & 0 & 0&\cdots& Q_{4s}^{t_1,t_2}& Q_{4s-1}^{t_1,t_2} \\ [4pt]
	0 & 0 & 0&\cdots& 0& Q_{4s}^{t_1,t_2}  \\
	\end{array}\right]
	$$}
where every $2$-form is orderly and diagonally displayed.  
Then the concatenation
$$A:=[B_{1,1}~B_{1,2}~B_{1,3}~B_{2,2}~B_{2,3}~B_{3,3}]$$
is a $(5s-1)\times 6(s-1)$ matrix of maximal rank.
\end{enumerate}
\end{Lemma}
\demo (a)
We claim that the $(5s+2)\times(5s+2)$ submatrix $B$ of $A$ omitting the first $(s-2)$ columns of block $B_{3,3}$ has nonzero determinant.

Note that the elements along the main diagonal of $B$ are
\begin{equation*}\nonumber
\underbrace{Q_{0}^{1,1}, \cdots, Q_{0}^{1,1}}_s,~\underbrace{Q_{s}^{1,2}, \cdots, Q_{s}^{1,2}}_s,~\underbrace{Q_{2s}^{1,3}, \cdots, Q_{2s}^{1,3}}_s,~\underbrace{Q_{3s}^{2,2}, \cdots, Q_{3s}^{2,2}}_s,~\underbrace{Q_{4s}^{2,3},\cdots , Q_{4s}^{2,3}}_s,~\underbrace{Q_{4s+2}^{3,3},~Q_{4s+2}^{3,3}}_2
\end{equation*}

We now specialize via the $k$-algebra endomorphism $\varphi$ of $S$ that fixes each one of the variables
$$T_0^{(1)}, T_{s}^{(2)}, T_{2s}^{(2)}, T_{2s+1}^{(3)}$$
and maps the remaining ones to $0$.

It suffices to show that the specialized matrix  has a nonzero determinant.

\pagebreak

The following table lists the only entries of $B$ that do not vanish under this specialization:

\bigskip

\begin{center}
\begin{tabular}{|c|c|c|}
\hline  The $2$-form & The relevant variables & the output\\
\hline $Q_0^{1,1}$ & $T_{0}^{(1)}, T_{0}^{(1)}$ & $\varphi(Q_0^{1,1})=T_{0}^{(1)}T_{0}^{(1)}=:q_0^{1,1}$    \\ 
\hline   $Q_s^{1,2}$ & $T_{0}^{(1)},T_{s}^{(2)}$ & $\varphi(Q_s^{1,2})=T_{0}^{(1)}T_{s}^{(2)}=:q_s^{1,2}$\\ 
\hline $Q_{2s}^{1,2}$ & $T_{0}^{(1)},T_{2s}^{(2)}$  &  $\varphi(Q_{2s}^{1,2})=T_{0}^{(1)}T_{2s}^{(2)}=:q_{2s}^{1,2}$\\  
\hline $Q_{2s+1}^{1,3}$ & $T_{0}^{(1)},T_{2s+1}^{(3)}$ & $\varphi(Q_{2s+1}^{1,3})=T_{0}^{(1)}T_{2s+1}^{(3)}=:q_{2s+1}^{1,3}$\\ 
\hline  $Q_{2s}^{2,2}$ & $T_{s}^{(2)},T_{s}^{(2)}$   & $\varphi(Q_{2s}^{2,2})=T_{s}^{(2)}T_{s}^{(2)}=:q_{2s}^{2,2}$\\ 
\hline $Q_{3s}^{2,2}$ & $T_{s}^{(2)},T_{2s}^{(2)}$   & $\varphi(Q_{3s}^{2,2})=2T_{s}^{(2)}T_{2s}^{(2)}=:q_{3s}^{2,2}$\\
\hline $Q_{3s+1}^{2,3}$ & $T_{s}^{(2)},T_{2s+1}^{(3)}$  & $\varphi(Q_{3s+1}^{2,3})=T_{s}^{(2)}T_{2s+1}^{(3)}=:q_{3s+1}^{2,3}$\\ 
\hline $Q_{4s}^{2,2}$ & $T_{2s}^{(2)},T_{2s}^{(2)}$   & $\varphi(Q_{4s}^{2,2})=T_{2s}^{(2)}T_{2s}^{(2)}=:q_{4s}^{2,2}$\\
\hline $Q_{4s+1}^{2,3}$ & $T_{2s}^{(2)},T_{2s+1}^{(3)}$ &  $\varphi(Q_{4s+1}^{2,3})=T_{2s}^{(2)}T_{2s+1}^{(3)}=:q_{4s+1}^{2,3}$\\ 
\hline $Q_{4s+2}^{3,3}$ & $T_{2s+1}^{(3)},T_{2s+1}^{(3)}$  &  $\varphi(Q_{4s+2}^{3,3})=T_{2s+1}^{(3)}T_{2s+1}^{(3)}=:q_{4s+2}^{3,3}$\\
\hline
\end{tabular}
\end{center}

\bigskip

Therefore, the resulting matrix $\widetilde{B}$ has the following shape, where now all nonzero entries are monomials of degree $2$ with coefficients $1$ or $2$:

$$\kern 0.3cm\overbrace{\hspace{1.7cm}}^{\widetilde{B}_{1,1}}\hspace{0.7cm}
\overbrace{\hspace{1.7cm}}^{\widetilde{B}_{1,2}}\hspace{0.7cm}
\overbrace{\hspace{2cm}}^{\widetilde{B}_{1,3}}\hspace{0.6cm}
\overbrace{\hspace{1.8cm}}^{\widetilde{B}_{2,2}}\hspace{0.6cm}
\overbrace{\hspace{1.8cm}}^{\widetilde{B}_{2,3}}\hspace{2.5cm} $$
$${\tiny\left[\begin{array}{ccccccccccccccccccc}
q_0^{1,1}  & \cdots & 0 & 0  & \cdots & 0& 0& \cdots& 0& 0 & \cdots & 0 & 0& \cdots& 0& 0& 0\\
\vdots   &  & \vdots & \vdots   & & \vdots&   \vdots &  &\vdots & \vdots  & & \vdots& \vdots& & \vdots& \vdots& \vdots\\
0  & \cdots & q_0^{1,1}& 0 &  \cdots& 0&  0&  \cdots& 0& 0&  \cdots& 0& 0& \cdots& 0& 0& 0\\
0  & \cdots & 0& q_{s}^{1,2} &  \cdots& 0&  0&  \cdots& 0& 0&  \cdots& 0& 0& \cdots& 0& 0& 0\\
\vdots   &  & \vdots & \vdots  &  & \vdots&   \vdots & &\vdots &\vdots   & & \vdots& \vdots& & \vdots& \vdots& \vdots\\
0  & \cdots & 0& 0 & \cdots& q_{s}^{1,2}&  0&  \cdots& 0& 0&  \cdots & 0& 0& \cdots& 0& 0& 0\\

0  & \cdots & 0& q_{2s}^{1,2} & \cdots& 0& 0 &  \cdots& 0& q_{2s}^{2,2}&  \cdots& 0& 0& \cdots& 0& 0& 0\\  [5pt]
0  & \cdots & 0& 0& \cdots& 0&  q_{2s+1}^{1,3}&  \cdots& 0& 0&  \cdots & 0& 0& \cdots& 0& 0& 0\\
\vdots   &  & \vdots & \vdots &   & \vdots  & \vdots & &\vdots  &\vdots  & & \vdots& \vdots& & \vdots& \vdots& \vdots\\ 

0  & \cdots & 0& 0 &  \cdots& q_{2s}^{1,2}&  0&  \cdots&  0& 0&  \cdots& q_{2s}^{2,2}& 0& \cdots& 0& 0& 0\\  [5pt]
0  & \cdots & 0& 0 &  \cdots& 0&  0&  \cdots&  q_{2s+1}^{1,3}& q_{3s}^{2,2}& \cdots& 0& 0& \cdots& 0& 0& 0\\
0  & \cdots & 0& 0 &  \cdots& 0&  0&  \cdots&  0& 0&  \cdots& 0& q_{3s+1}^{2,3}& \cdots& 0& 0& 0\\
0  & \cdots & 0& 0 &  \cdots& 0&  0&  \cdots&  0& 0&  \cdots& 0& 0& \cdots& 0& 0& 0\\
\vdots   &  & \vdots & \vdots  &  & \vdots  & \vdots & &\vdots  &\vdots  & & \vdots& \vdots& & \vdots& \vdots& \vdots\\
0  & \cdots & 0& 0 &  \cdots& 0&  0&  \cdots&  0& 0&  \cdots& q_{3s}^{2,2}& 0& \cdots& 0& 0& 0\\
0  & \cdots & 0& 0 &  \cdots& 0&  0&  \cdots&  0& q_{4s}^{2,2}&  \cdots& 0& 0& \cdots& q_{3s+1}^{2,3}& 0& 0\\ [5pt]
0  & \cdots & 0& 0 &  \cdots& 0&  0&  \cdots&  0& 0&  \cdots& 0& q_{4s+1}^{2,3}&\cdots& 0& 0& 0\\
\vdots   &  & \vdots & \vdots  &  & \vdots  & \vdots & &\vdots  &\vdots & & \vdots& \vdots& & \vdots& \vdots& \vdots\\

0  & \cdots & 0& 0 &  \cdots& 0&  0&  \cdots&  0& 0& \cdots & q_{4s}^{2,2}& 0& \cdots& 0& 0& 0\\  [5pt]
0  & \cdots & 0& 0 &  \cdots& 0&  0&  \cdots&  0& 0&  \cdots & 0& 0& \cdots& q_{4s+1}^{2,3}& q_{4s+2}^{3,3}& 0\\  [5pt]
0  & \cdots & 0& 0 &  \cdots& 0&  0&  \cdots&  0& 0&  \cdots & 0& 0& \cdots& 0& 0&q_{4s+2}^{3,3}
\end{array}\right]} $$

\bigskip
Accordingly, the elements along the main diagonal of the above matrix are

\bigskip
\begin{equation*}\nonumber
\underbrace{q_{0}^{1,1}, \cdots, q_{0}^{1,1}}_s,~\underbrace{q_{s}^{1,2}, \cdots, q_{s}^{1,2}}_s,~\underbrace{0, \cdots, 0}_s,~\underbrace{q_{3s}^{2,2}, \cdots, q_{3s}^{2,2}}_s,~\underbrace{0,\cdots , 0}_s,~\underbrace{q_{4s+2}^{3,3},~q_{4s+2}^{3,3}}_2
\end{equation*}

In order to prove that the determinant is nonzero we use a strategy of blocks.
To start,  all the entries to the right of the slots of either one of the entries $q_{0}^{1,1}$ and $q_{s}^{1,2}$ are null and, likewise, so are the entries to the left of the rightmost slot of $q_{4s+2}^{3,3}$.
Therefore, we are reduced to showing the non vanishing of the determinant of the following matrix

$${\small\left[\begin{array}{ccccccccccccccccccc}

 0 &  \cdots& 0& q_{2s}^{2,2}&  \cdots& 0& 0& \cdots& 0&0& 0\\  [5pt]
 q_{2s+1}^{1,3}&  \cdots& 0& 0&  \cdots & 0& 0& \cdots&0& 0& 0\\
 \vdots & &\vdots  &\vdots  & & \vdots& \vdots& &\vdots& \vdots& \vdots \\ 

  0&  \cdots&  0& 0&  \cdots& q_{2s}^{2,2}& 0& \cdots& 0&0& 0\\  [5pt]
  0&  \cdots&  q_{2s+1}^{1,3}& q_{3s}^{2,2}& \cdots& 0& 0& \cdots& 0&0& 0\\
  0&  \cdots&  0& 0&  \cdots& 0& q_{3s+1}^{2,3}& \cdots& 0&0& 0\\
  0&  \cdots&  0& 0&  \cdots& 0& 0& \cdots& 0&0& 0\\
 \vdots & &\vdots  &\vdots  & & \vdots& \vdots& & \vdots&\vdots& \vdots \\
  0&  \cdots&  0& 0&  \cdots& q_{3s}^{2,2}& 0& \cdots& q_{3s+1}^{2,3}&0& 0\\
  0&  \cdots&  0& q_{4s}^{2,2}&  \cdots& 0& 0& \cdots& 0&q_{3s+1}^{2,3}& 0\\ [5pt]
  0&  \cdots&  0& 0&  \cdots& 0& q_{4s+1}^{2,3}&\cdots& 0&0& 0\\
 \vdots & &\vdots  &\vdots & & \vdots& \vdots& &\vdots& \vdots& \vdots\\

  0&  \cdots&  0& 0& \cdots & q_{4s}^{2,2}& 0& \cdots& q_{4s+1}^{2,3}& 0& 0\\  [5pt]
  0&  \cdots&  0& 0&  \cdots & 0& 0& \cdots& 0& q_{4s+1}^{2,3}& q_{4s+2}^{3,3}\\  [5pt]
\end{array}\right]} $$

Similarly for this matrix, the entries above the entry $q_{4s+2}^{3,3}$ are all zero, hence it suffices to show that the following submatrix has a nonzero determinant

$${\small\left[\begin{array}{ccccccccccccccccccc}

 0 &  \cdots& 0& q_{2s}^{2,2}& 0& \cdots& 0& 0& \cdots& 0& 0\\  [5pt]
 q_{2s+1}^{1,3}&  \cdots& 0& 0& q_{2s}^{2,2}& \cdots & 0& 0& \cdots& 0& 0\\
 \vdots & &\vdots  &\vdots &\vdots & & \vdots& \vdots& &\vdots& \vdots\\ 

  0&  \cdots&  0& 0& 0& \cdots& q_{2s}^{2,2}& 0& \cdots& 0& 0\\  [5pt]
  0&  \cdots&  q_{2s+1}^{1,3}& q_{3s}^{2,2}& 0&\cdots& 0& 0& \cdots& 0& 0\\
  0&  \cdots&  0& 0& q_{3s}^{2,2}& \cdots& 0& q_{3s+1}^{2,3}& \cdots& 0& 0\\
 
 \vdots & &\vdots  &\vdots &\vdots  & & \vdots& \vdots& &\vdots& \vdots \\
  0&  \cdots&  0& 0& 0&  \cdots& q_{3s}^{2,2}& 0& \cdots& q_{3s+1}^{2,3}& 0\\
  0&  \cdots&  0& q_{4s}^{2,2}&0&  \cdots& 0& 0& \cdots& 0& q_{3s+1}^{2,3}\\ [5pt]
  0&  \cdots&  0& 0& q_{4s}^{2,2}& \cdots& 0& q_{4s+1}^{2,3}&\cdots& 0& 0\\
 \vdots & &\vdots  &\vdots &\vdots & & \vdots& \vdots& &\vdots& \vdots\\

  0&  \cdots&  0& 0& 0&\cdots & q_{4s}^{2,2}& 0& \cdots& q_{4s+1}^{2,3}& 0\\  [5pt]
\end{array}\right]} $$

Next observe that the entries above and under each slot of $q^{1,3}_{2s+1}$  are null. Then, we take these slots successively as pivots along their common diagonal; and subsequently, the leftmost slot of $q_{2s}^{2,2}$ as pivot, thus reducing to the matrix

$${\normalsize\left[\begin{array}{ccccccccccccccccccc}

  q_{3s}^{2,2}&\cdots& 0& q_{3s+1}^{2,3}& \cdots& 0& 0\\
 
\vdots & & \vdots& \vdots& &\vdots& \vdots \\
  0&\cdots& q_{3s}^{2,2}& 0& \cdots& q_{3s+1}^{2,3}& 0\\
  0&\cdots& 0& 0& \cdots& 0& q_{3s+1}^{2,3}\\ [5pt]
 q_{4s}^{2,2} &\cdots& 0& q_{4s+1}^{2,3}&\cdots& 0& 0\\
\vdots & & \vdots& \vdots& &\vdots& \vdots\\

 0& \cdots & q_{4s}^{2,2}& 0& \cdots& q_{4s+1}^{2,3}& 0\\  [5pt]

\end{array}\right]} $$

Now, moving the row having $q^{2,3}_{3s+1}$ as rightmost entry all the way down to the last row yields the matrix

$${\normalsize\left[\begin{array}{ccccccccccccccccccc}

  q_{3s}^{2,2}&\cdots& 0& q_{3s+1}^{2,3}& \cdots& 0& 0\\
 
\vdots & & \vdots& \vdots& &\vdots& \vdots \\
  0&\cdots& q_{3s}^{2,2}& 0& \cdots& q_{3s+1}^{2,3}& 0\\

 q_{4s}^{2,2} &\cdots& 0& q_{4s+1}^{2,3}&\cdots& 0& 0\\
\vdots & & \vdots& \vdots& &\vdots& \vdots\\

 0& \cdots & q_{4s}^{2,2}& 0& \cdots& q_{4s+1}^{2,3}& 0\\  [5pt]
  0&\cdots& 0& 0& \cdots& 0& q_{3s+1}^{2,3}\\ [5pt]
\end{array}\right]} $$

Finally, note that over the field of fractions of $S$ this matrix is equivalent to an upper triangular matrix with nonzero entries along the main diagonal.
Indeed, for this it suffices to use row operations depending only on the following non-vanishing $2\times 2$ determinant
$$q_{4s}^{2,2}q_{3s+1}^{2,3}-q_{3s}^{2,2}q_{4s+1}^{2,3}=T_{2s}^{(2)}T_{2s}^{(2)}T_{s}^{(2)}T_{2s+1}^{(3)}-2T_{s}^{(2)}T_{2s}^{(2)}T_{2s}^{(2)}T_{2s+1}^{(3)}=-T_{s}^{(2)}T_{2s}^{(2)}T_{2s}^{(2)}T_{2s+1}^{(3)}\neq 0.$$

\smallskip

(b) We claim that the $(5s-1)\times(5s-1)$ submatrix $B$ of $A$ omitting the first $(s-5)$ columns of block $B_{3,3}$ has nonzero determinant.

Note that the elements along the main diagonal of $B$ are

\begin{eqnarray*} \nonumber
\underbrace{Q_{0}^{1,1}, \cdots, Q_{0}^{1,1}}_{s-1},~\underbrace{Q_{s-1}^{1,2}, \cdots, Q_{s-1}^{1,2}}_{s-1},~\underbrace{Q_{2s-2}^{1,3}, \cdots, Q_{2s-2}^{1,3}}_{s-1},~\underbrace{Q_{3s-3}^{2,2}, \cdots, Q_{3s-3}^{2,2}}_{s-1},  \\ \nonumber
\underbrace{Q_{4s-4}^{2,3},  \cdots , Q_{4s-4}^{2,3}}_{s-1},~\underbrace{Q_{4s}^{3,3},~Q_{4s}^{3,3},~Q_{4s}^{3,3},~Q_{4s}^{3,3}}_4
\end{eqnarray*}

We now specialize via the $k$-algebra endomorphism $\varphi$ of $S$ that fixes each one of the variables
$$T_0^{(1)}, T_{s-1}^{(2)}, T_{2(s-1)}^{(2)}, T_{2s}^{(3)}$$
and maps the remaining ones to $0$.
As before, it suffices to show that the specialized matrix $\widetilde{B}$ has a nonzero determinant.

The following table lists the entries of $B$ that do not vanish under this specialization:

\medskip

\begin{center}
\begin{tabular}{|c|c|c|}
	\hline  The $2$-form & The relevant variables & the output\\
\hline  $Q_0^{1,1}$ &  $T_{0}^{(1)},T_{0}^{(1)}$ & $\varphi(Q_0^{1,1})=T_{0}^{(1)}T_{0}^{(1)}=:q_0^{1,1}$    \\ 
\hline $Q_{s-1}^{1,2}$ & $T_{0}^{(1)},T_{s-1}^{(2)}$ &   $\varphi(Q_{s-1}^{1,2})=T_{0}^{(1)}T_{s-1}^{(2)}=:q_{s-1}^{1,2}$\\ 
\hline $Q_{2(s-1)}^{1,2}$ & $T_{0}^{(1)},T_{2(s-1)}^{(2)}$ &  $\varphi(Q_{2(s-1)}^{1,2})=T_{0}^{(1)}T_{2(s-1)}^{(2)}=:q_{2(s-1)}^{1,2}$\\ 
\hline  $Q_{2s}^{1,3}$ & $T_{0}^{(1)},T_{2s}^{(3)}$ &  $\varphi(Q_{2s}^{1,3})=T_{0}^{(1)}T_{2s}^{(3)}=:q_{2s}^{1,3}$\\ 
\hline $Q_{2(s-1)}^{2,2}$ & $T_{s-1}^{(2)},T_{s-1}^{(2)}$ &   $\varphi(Q_{2(s-1)}^{2,2})=T_{s-1}^{(2)}T_{s-1}^{(2)}=:q_{2(s-1)}^{2,2}$\\ 
\hline  $Q_{3(s-1)}^{2,2}$& $T_{s-1}^{(2)},T_{2(s-1)}^{(2)}$ &  $\varphi(Q_{3(s-1)}^{2,2})=2T_{s-1}^{(2)}T_{2(s-1)}^{(2)}=:q_{3(s-1)}^{2,2}$\\
\hline $Q_{3s-1}^{2,3}$ & $T_{s-1}^{(2)},T_{2s}^{(3)}$ &  $\varphi(Q_{3s-1}^{2,3})=T_{s-1}^{(2)}T_{2s}^{(3)}=:q_{3s-1}^{2,3}$\\ 
\hline $Q_{4(s-1)}^{2,2}$ &$T_{2(s-1)}^{(2)}, T_{2(s-1)}^{(2)}$ &  $\varphi(Q_{4(s-1)}^{2,2})=T_{2(s-1)}^{(2)}T_{2(s-1)}^{(2)}=:q_{4(s-1)}^{2,2}$\\ 
\hline  $Q_{4s-2}^{2,3}$ & $T_{2(s-1)}^{(2)},T_{2s}^{(3)}$ &  $\varphi(Q_{4s-2}^{2,3})=T_{2(s-1)}^{(2)}T_{2s}^{(3)}=:q_{4s-2}^{2,3}$\\ 
\hline  $Q_{4s}^{3,3}$ & $T_{2s}^{(3)},T_{2s}^{(3)}$ &  $\varphi(Q_{4s}^{3,3})=T_{2s}^{(3)}T_{2s}^{(3)}=:q_{4s}^{3,3}$\\
\hline
\end{tabular}
\end{center}

Therefore,  $\widetilde{B}$ has the following shape, where again the nonzero entries are squarefree monomials of degree $2$ with coefficients $1$ or $2$.

\begin{center}
\begin{figure}[!htb]
\centering
\includegraphics[scale=0.76]{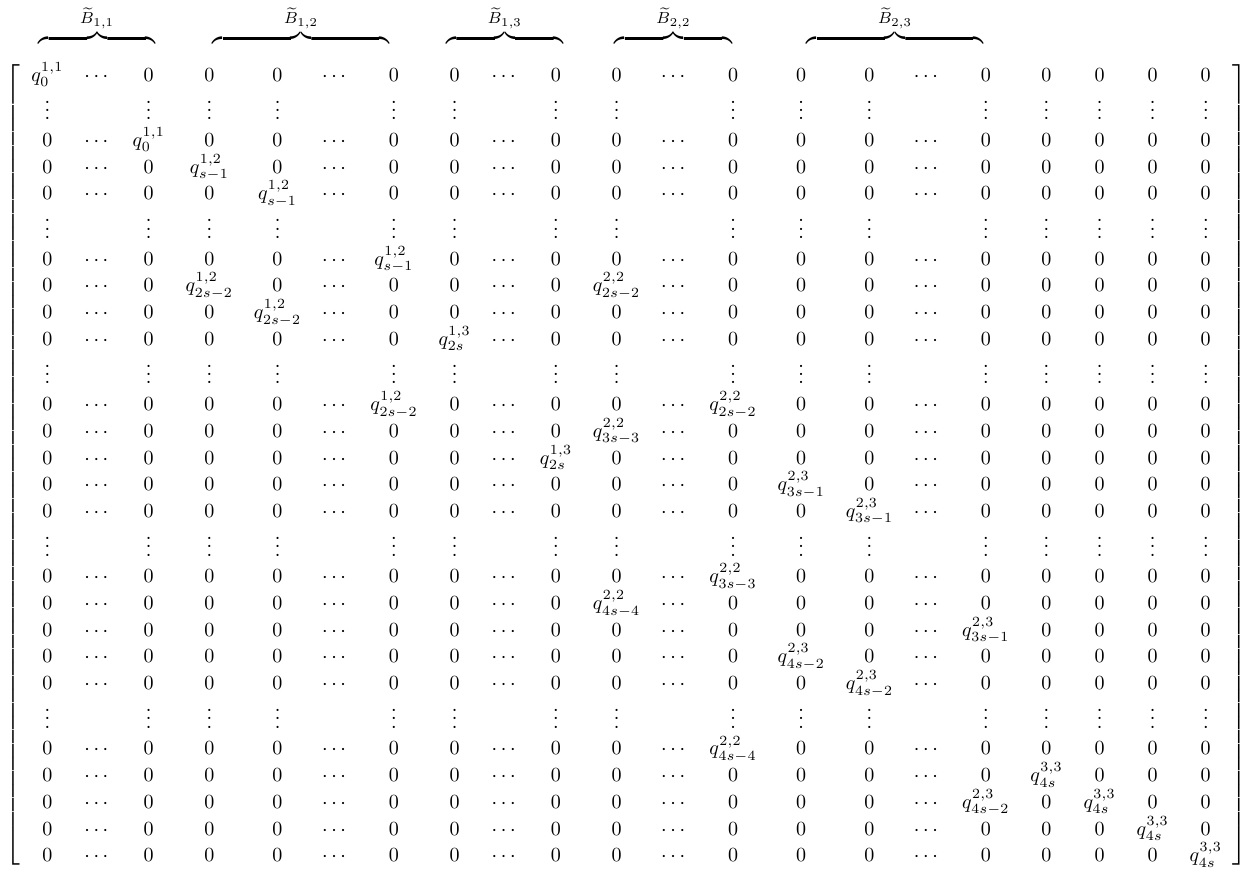}
\end{figure}
\end{center}
\pagebreak

Note that the diagonal elements of the submatrix above are

\begin{eqnarray*} \nonumber
\underbrace{q_{0}^{1,1}, \cdots, q_{0}^{1,1}}_{s-1},~\underbrace{q_{s-1}^{1,2}, \cdots, q_{s-1}^{1,2}}_{s-1},~\underbrace{0, \cdots, 0}_{s-1},~\underbrace{q_{3s-3}^{2,2}, \cdots, q_{3s-3}^{2,2}}_{s-1},
\underbrace{0,  \cdots , 0}_{s-1},~\underbrace{q_{4s}^{3,3},~q_{4s}^{3,3},~q_{4s}^{3,3},~q_{4s}^{3,3}}_4
\end{eqnarray*}
One proceeds in a completely analogous to the odd case,  using a strategy of blocks.
Often, as in the above matrix and in some of the subsequent ones, in order to improve visualization we choose to show more columns in certain blocks as compared to others.

To start,  all entries to the right of the slots of either one of the entries $q_{0}^{1,1}$ and $q_{s-1}^{1,2}$ are null and, likewise, so are all entries to the left of two rightmost slot of $q_{4s}^{3,3}$.
Therefore, we are reduced to showing the non vanishing of the determinant of the following matrix:

$${\tiny\left[\begin{array}{ccccccccccccccccccccc}
 0 &  \cdots& 0& q_{2s-2}^{2,2}& 0& \cdots& 0& 0& 0&\cdots& 0 & 0 & 0& 0& 0 \\
  0&  \cdots& 0& 0& q_{2s-2}^{2,2} & \cdots   & 0& 0& 0&\cdots& 0& 0& 0& 0& 0 \\
  q_{2s}^{1,3}&  \cdots& 0& 0&  0&\cdots & 0& 0& 0&\cdots& 0& 0& 0 & 0& 0\\

 \vdots & &\vdots  &\vdots  & \vdots& & \vdots& \vdots& \vdots&& \vdots& \vdots& \vdots& \vdots& \vdots\\

  0&  \cdots&  0& 0&  0& \cdots& q_{2s-2}^{2,2}& 0& 0&\cdots& 0& 0& 0 & 0& 0\\
  0&  \cdots&  0 & q_{3s-3}^{2,2}& 0& \cdots& 0& 0& 0&\cdots& 0& 0& 0& 0& 0\\
  0&  \cdots&  q_{2s}^{1,3}& 0& q_{3s-3}^{2,2} & \cdots& 0& 0& 0&\cdots& 0& 0& 0& 0& 0\\
  0&  \cdots&  0& 0&  0& \cdots& 0& q_{3s-1}^{2,3}& 0&\cdots& 0& 0& 0& 0& 0\\
  0&  \cdots&  0& 0&  0& \cdots& 0& 0& q_{3s-1}^{2,3}&\cdots& 0& 0& 0& 0& 0\\

\vdots & &\vdots  &\vdots &\vdots  & & \vdots& \vdots& \vdots&& \vdots& \vdots& \vdots& \vdots& \vdots\\

  0&  \cdots&  0& 0& 0& \cdots& q_{3s-3}^{2,2}& 0& 0&\cdots& q_{3s-1}^{2,3}& 0& 0& 0& 0\\
  0&  \cdots&  0& q_{4s-4}^{2,2}& 0&  \cdots& 0& 0& 0&\cdots& 0& q_{3s-1}^{2,3}& 0& 0& 0\\
  0&  \cdots&  0& 0&  q_{4s-4}^{2,2}& \cdots& 0& 0& 0&\cdots& 0& 0& q_{3s-1}^{2,3}& 0& 0\\
  0&  \cdots&  0& 0& 0& \cdots& 0& q_{4s-2}^{2,3}&0&\cdots& 0& 0& 0& 0& 0\\
  0&  \cdots&  0& 0& 0& \cdots& 0& 0&q_{4s-2}^{2,3}&\cdots& 0& 0& 0& 0& 0\\

 \vdots & &\vdots  &\vdots &\vdots & & \vdots& \vdots&\vdots& & \vdots& \vdots& \vdots& \vdots& \vdots\\

  0&  \cdots&  0& 0& 0&\cdots & q_{4s-4}^{2,2}& 0& 0&\cdots& q_{4s-2}^{2,3}& 0& 0& 0& 0\\
  0&  \cdots&  0& 0& 0&\cdots & 0& 0& 0&\cdots& 0& q_{4s-2}^{2,3}& 0& q_{4s}^{3,3}& 0\\
  0&  \cdots&  0& 0&  0&\cdots & 0& 0& 0&\cdots& 0&  0& q_{4s-2}^{2,3} & 0& q_{4s}^{3,3}\\
\end{array}\right]} $$

And again, the entries above the two slots of the entriy $q_{4s}^{3,3}$ are all zero, hence it suffices to show that the following submatrix has a nonzero determinant:

$${\tiny\left[\begin{array}{ccccccccccccccccccccc}

 0 &  \cdots& 0& q_{2s-2}^{2,2}& 0& 0& \cdots& 0& 0& 0&\cdots& 0 & 0 & 0 & 0  \\
  0&  \cdots& 0& 0& q_{2s-2}^{2,2} & 0 &\cdots   & 0& 0& 0&\cdots& 0& 0& 0 & 0 \\
  q_{2s}^{1,3}&  \cdots& 0& 0&  0& q_{2s-2}^{2,2} &\cdots & 0& 0& 0&\cdots& 0& 0& 0 & 0 \\

 \vdots & &\vdots  &\vdots  & \vdots &\vdots& & \vdots& \vdots& \vdots&& \vdots& \vdots& \vdots& \vdots\\

  0&  \cdots&  0& 0&  0& 0&\cdots& q_{2s-2}^{2,2}& 0& 0&\cdots& 0& 0& 0 & 0 \\
  0&  \cdots&  0 & q_{3s-3}^{2,2}& 0& 0& \cdots& 0& 0& 0&\cdots& 0& 0& 0& 0 \\
  0&  \cdots&  q_{2s}^{1,3}& 0& q_{3s-3}^{2,2} & 0& \cdots& 0& 0& 0&\cdots& 0& 0& 0& 0 \\
  0&  \cdots&  0& 0&  0& q_{3s-3}^{2,2} & \cdots& 0& q_{3s-1}^{2,3}& 0&\cdots& 0& 0& 0& 0 \\
  0&  \cdots&  0& 0&  0& 0&\cdots& 0& 0& q_{3s-1}^{2,3}&\cdots& 0& 0& 0& 0 \\

\vdots & &\vdots  &\vdots &\vdots &\vdots  & & \vdots& \vdots& \vdots&& \vdots& \vdots& \vdots& \vdots\\

  0&  \cdots&  0& 0& 0& 0& \cdots& 0& 0& 0&\cdots& q_{3s-1}^{2,3}& 0& 0& 0 \\
  0&  \cdots&  0& 0& 0& 0& \cdots& q_{3s-3}^{2,2}& 0& 0&\cdots& 0& q_{3s-1}^{2,3}& 0& 0 \\
  0&  \cdots&  0& q_{4s-4}^{2,2}& 0& 0&  \cdots& 0& 0& 0&\cdots& 0 & 0& q_{3s-1}^{2,3}& 0\\
  0&  \cdots&  0& 0&  q_{4s-4}^{2,2}& 0 &\cdots& 0& 0& 0&\cdots& 0& 0& 0& q_{3s-1}^{2,3}\\
  0&  \cdots&  0& 0& 0& q_{4s-4}^{2,2}&\cdots& 0& q_{4s-2}^{2,3}&0&\cdots& 0& 0& 0& 0\\
  0&  \cdots&  0& 0& 0& 0&\cdots& 0& 0& q_{4s-2}^{2,3}&\cdots& 0& 0&0& 0\\

 \vdots & &\vdots  &\vdots &\vdots &\vdots & & \vdots& \vdots&\vdots& & \vdots&  \vdots &\vdots& \vdots\\

  0&  \cdots&  0& 0& 0& 0 &\cdots & 0& 0& 0&\cdots& q_{4s-2}^{2,3}& 0& 0& 0\\
  0&  \cdots&  0& 0& 0& 0 &\cdots & q_{4s-4}^{2,2}& 0& 0&\cdots&0& q_{4s-2}^{2,3}& 0& 0\\

\end{array}\right]} $$

Proceeding in this way, now take as pivot all slots of $q^{1,3}_{2s}$, noting that the entries above and under these slots are null, and subsequently, the two leftmost slots of $q_{2s-2}^{2,2}$, thus reducing to the matrix

$${\normalsize\left[\begin{array}{ccccccccccccccccccccc}

 q_{3s-3}^{2,2} & \cdots& 0& q_{3s-1}^{2,3}& 0&\cdots& 0& 0& 0& 0 \\
 0&\cdots& 0& 0& q_{3s-1}^{2,3}&\cdots& 0& 0& 0& 0 \\

\vdots  & & \vdots& \vdots& \vdots&& \vdots& \vdots& \vdots& \vdots\\

 0& \cdots& 0& 0& 0&\cdots& q_{3s-1}^{2,3}& 0& 0& 0 \\
 0& \cdots& q_{3s-3}^{2,2}& 0& 0&\cdots& 0& q_{3s-1}^{2,3}& 0& 0 \\
 0&  \cdots& 0& 0& 0&\cdots& 0 & 0& q_{3s-1}^{2,3}& 0\\
 0 &\cdots& 0& 0& 0&\cdots& 0& 0& 0& q_{3s-1}^{2,3}\\
 q_{4s-4}^{2,2}&\cdots& 0& q_{4s-2}^{2,3}&0&\cdots& 0& 0& 0& 0\\
 0&\cdots& 0& 0& q_{4s-2}^{2,3}&\cdots& 0& 0&0& 0\\

\vdots & & \vdots& \vdots&\vdots& & \vdots&  \vdots &\vdots& \vdots\\

 0 &\cdots & 0& 0& 0&\cdots& q_{4s-2}^{2,3}& 0& 0& 0\\
 0 &\cdots & q_{4s-4}^{2,2}& 0& 0&\cdots&0& q_{4s-2}^{2,3}& 0& 0\\

\end{array}\right]} $$

Now, moving the two row having $q^{2,3}_{3s-1}$ as rightmost entry all the way down to the last two rows yields the matrix

$${\normalsize\left[\begin{array}{ccccccccccccccccccccc}

 q_{3s-3}^{2,2} & \cdots& 0& q_{3s-1}^{2,3}& 0&\cdots& 0& 0& 0& 0 \\
 0&\cdots& 0& 0& q_{3s-1}^{2,3}&\cdots& 0& 0& 0& 0 \\

\vdots  & & \vdots& \vdots& \vdots&& \vdots& \vdots& \vdots& \vdots\\

 0& \cdots& 0& 0& 0&\cdots& q_{3s-1}^{2,3}& 0& 0& 0 \\
 0& \cdots& q_{3s-3}^{2,2}& 0& 0&\cdots& 0& q_{3s-1}^{2,3}& 0& 0 \\

 q_{4s-4}^{2,2}&\cdots& 0& q_{4s-2}^{2,3}&0&\cdots& 0& 0& 0& 0\\
 0&\cdots& 0& 0& q_{4s-2}^{2,3}&\cdots& 0& 0&0& 0\\

\vdots & & \vdots& \vdots&\vdots& & \vdots&  \vdots &\vdots& \vdots\\

 0 &\cdots & 0& 0& 0&\cdots& q_{4s-2}^{2,3}& 0& 0& 0\\
 0 &\cdots & q_{4s-4}^{2,2}& 0& 0&\cdots&0& q_{4s-2}^{2,3}& 0& 0\\

 0&  \cdots& 0& 0& 0&\cdots& 0 & 0& q_{3s-1}^{2,3}& 0\\
 0 &\cdots& 0& 0& 0&\cdots& 0& 0& 0& q_{3s-1}^{2,3}\\

\end{array}\right]} $$

Finally, note that over the field of fractions of $S$ this matrix is equivalent to an upper triangular matrix with nonzero entries along the main diagonal.
Indeed, for this it suffices to use row operations depending only on the following non-vanishing $2\times 2$ determinant
\begin{eqnarray}\nonumber
q_{4s-4}^{2,2}q_{3s-1}^{2,3}-q_{3s-3}^{2,2}q_{4s-2}^{2,3}&=&T_{2s-2}^{(2)}T_{2s-2}^{(2)}T_{s-1}^{(2)}T_{2s}^{(3)}-2T_{s-1}^{(2)}T_{2s-2}^{(2)}T_{2s-2}^{(2)}T_{2s}^{(3)}\\ \nonumber
&=&-T_{s-1}^{(2)}T_{2s-2}^{(2)}T_{2s-2}^{(2)}T_{2s}^{(3)}\neq 0 
\end{eqnarray}
\qed

\section{The main argument}

In this section we examine a certain crucial annihilator related to the Ratliff--Rush filtrations of $I$ and its powers -- the latter are known to give a bound for the depth of the associated graded ring of $I$.

The method employed here hinges on Lemma~\ref{main_lemma} and is entirely explicit.

\begin{Proposition}\label{annihilators} Let $I\subset R:=k[x,y]$ denote a codimension $2$ ideal generated by a set of general forms $f_1,f_2,f_3$ of degree $d\geq 5$ {\rm (}respectively, $d\geq 10${\rm )} if $d$ is odd {\rm (}respectively, even{\rm )}.
Then $I^2:I\not\subset I $.
\end{Proposition}
\demo
Since a change of variables does not affect either the hypothesis or the statement, we can assume by Lemma~\ref{Hilbert_function}  that

$\bullet$ If $d=2s+1$, then $y^{3s}$ spans $(R/I)_{3s}$; 

$\bullet$  If $d=2s$, then $\{y^{3s-2}, xy^{3s-3}\}$ span $(R/I)_{3s-2}$.

 









\noindent {\bf Odd case:} $d=2s+1$. 

We are assuming that $s\geq 2$.
We are to show that $y^{3s}f_j\in I^2$, for $j=1,2,3$.

\smallskip

Writing $g_{1,1}=(f_1)^2,~g_{1,2}=f_1f_2,~g_{1,3}=f_1f_3,~g_{2,2}=(f_2)^2,~g_{2,3}=f_2f_3,~g_{3,3}=(f_3)^2$, we are asking for a solution of the equation 
 \begin{equation}\label{content}
y^{3s}f_1=\boldsymbol{h}_{1,1}g_{1,1}+\cdots+\boldsymbol{h}_{3,3}g_{3,3},
\end{equation} 
in homogeneous forms  $\boldsymbol{h}_{t_1,t_2}$ of degree $s-1$, $1\leq t_1\leq t_2\leq 3$.

\smallskip

Set
\begin{eqnarray*}\nonumber
f_j&=&\gamma_{0}^{(j)}y^{2s+1}+\cdots+\gamma_{2s+1}^{(j)}x^{2s+1} \\ \nonumber
g_{t_1,t_2}&=&\beta_{0}^{t_1,t_2}y^{4s+2}+\cdots+\beta_{4s+2}^{t_1,t_2}x^{4s+2}\\ \nonumber
\boldsymbol{h}_{t_1,t_2}&=&\alpha_{0}^{t_1,t_2}y^{s-1}+\cdots+\alpha_{s-1}^{t_1,t_2}x^{s-1}
\end{eqnarray*}
where $1\leq j\leq 3$ and $1\leq t_1\leq t_2\leq 3$. Here $\{\gamma_{l}^{(j)}\}$ are general coefficients in the ground field $k$ and $\{\beta_{l'}^{t_1,t_2}\}$ are the resulting quadratic expressions, while $\{\alpha_i\,|\, 0\leq i \leq s-1\}$ are sought for.

From the relation in (\ref{content}), comparing equal monomials in $x,y$ on both sides, the following $5s+2$ relations are obtained

$$\begin{array}{cccccc}
&\gamma_{0}^{(1)} & = & \displaystyle\sum_{(t_1,t_2)=(1,1)}^{(3,3)}\alpha_{0}^{t_1,t_2}\beta_{0}^{t_1,t_2} \hspace{10cm}\\ 
&\gamma_{1}^{(1)} & = & \displaystyle\sum_{(t_1,t_2)=(1,1)}^{(3,3)}\alpha_{0}^{t_1,t_2}\beta_{1}^{t_1,t_2}+\alpha_{1}^{t_1,t_2}\beta_{0}^{t_1,t_2}\hspace{7.8cm} \\ 
&\gamma_{2}^{(1)} & = & \displaystyle\sum_{(t_1,t_2)=(1,1)}^{(3,3)}\alpha_{0}^{t_1,t_2}\beta_{2}^{t_1,t_2}+\alpha_{1}^{t_1,t_2}\beta_{1}^{t_1,t_2}+\alpha_{2}^{t_1,t_2}\beta_{0}^{t_1,t_2} \hspace{5.6cm}\\
&\gamma_{3}^{(1)} & = & \displaystyle\sum_{(t_1,t_2)=(1,1)}^{(3,3)}\alpha_{0}^{t_1,t_2}\beta_{3}^{t_1,t_2}+\alpha_{1}^{t_1,t_2}\beta_{2}^{t_1,t_2}+\alpha_{2}^{t_1,t_2}\beta_{1}^{t_1,t_2}+\alpha_{3}^{t_1,t_2}\beta_{0}^{t_1,t_2} \hspace{3.4cm}\\
&\vdots  && \vdots \hspace{7cm}\\
&\gamma_{s-1}^{(1)} & = &  \displaystyle\sum_{(t_1,t_2)=(1,1)}^{(3,3)}\alpha_{0}^{t_1,t_2}\beta_{s-1}^{t_1,t_2}+\cdots+\alpha_{s-1}^{t_1,t_2}\beta_{0}^{t_1,t_2}  \hspace{6.8cm}\\ 
&\gamma_{s}^{(1)} & = &  \displaystyle\sum_{(t_1,t_2)=(1,1)}^{(3,3)}\alpha_{0}^{t_1,t_2}\beta_{s}^{t_1,t_2}+\cdots+\alpha_{s-1}^{t_1,t_2}\beta_{1}^{t_1,t_2}\hspace{6.8cm}\\ 
&\vdots  & & \vdots \hspace{7cm}\\
&\gamma_{2s+1}^{(1)}& = &  \displaystyle\sum_{(t_1,t_2)=(1,1)}^{(3,3)}\alpha_{0}^{t_1,t_2}\beta_{2s+1}^{t_1,t_2}+\cdots+\alpha_{s-1}^{t_1,t_2}\beta_{s+2}^{t_1,t_2}\hspace{6.8cm}\\
\end{array}$$ 

$$
\begin{array}{cccc} 
&0 & = &  \displaystyle\sum_{(t_1,t_2)=(1,1)}^{(3,3)}\alpha_{0}^{t_1,t_2}\beta_{2s+2}^{t_1,t_2}+\cdots+\alpha_{s-1}^{t_1,t_2}\beta_{s+3}^{t_1,t_2}\hspace{6.8cm}\\ 
&\vdots &  &  \hspace{7cm}\\
&0 & = &  \displaystyle\sum_{(t_1,t_2)=(1,1)}^{(3,3)}\alpha_{0}^{t_1,t_2}\beta_{4s+1}^{t_1,t_2}+\cdots+\alpha_{s-1}^{t_1,t_2}\beta_{3s+2}^{t_1,t_2}\hspace{6.8cm}\\ 
&0& = &  \displaystyle\sum_{(t_1,t_2)=(1,1)}^{(3,3)}\alpha_{0}^{t_1,t_2}\beta_{4s+2}^{t_1,t_2}+\cdots+\alpha_{s-1}^{t_1,t_2}\beta_{3s+3}^{t_1,t_2}\hspace{6.8cm}\\ 
&0 & = &  \displaystyle\sum_{(t_1,t_2)=(1,1)}^{(3,3)}\alpha_{1}^{t_1,t_2}\beta_{4s+2}^{t_1,t_2}+\cdots+\alpha_{s-1}^{t_1,t_2}\beta_{3s+4}^{t_1,t_2}\hspace{6.8cm}\\
&\vdots &  &  \hspace{7cm}\\
&0 & = & \displaystyle\sum_{(t_1,t_2)=(1,1)}^{(3,3)}\alpha_{s-1}^{t_1,t_2}\beta_{4s+2}^{t_1,t_2}\hspace{9.8cm}
\end{array},
$$ 
where the vanishing of the left side member happens for the last $5s+1-(2s+2)+1=3s$  relations because $f_1$ has degree $2s+1.$ 

\smallskip

The  $(5s+2)\times 1$ vector on the left side is given and we look for a solution in the $\alpha$'s of the linear system defined by the following $(5s+2)\times 6s$ matrix

$$\small \mathcal A=\left[\begin{array}{cccccccccccccccccccc}
\beta_{0}^{1,1}& 0& 0& \cdots& 0 & 0&  &\cdots\cdots &  & \beta_{0}^{3,3} & 0& 0& \cdots & 0 & 0  \\ [5pt]
\beta_{1}^{1,1}& \beta_{0}^{1,1}& 0& \cdots& 0 &0&   &\cdots\cdots &   & \beta_{1}^{3,3} & \beta_{0}^{3,3}& 0& \cdots & 0 & 0  \\ [5pt]
\beta_{2}^{1,1}&\beta_{1}^{1,1}&\beta_{0}^{1,1}&\cdots & 0& 0&&\cdots\cdots&&\beta_{2}^{3,3}&\beta_{1}^{3,3}&\beta_{0}^{3,3}&\cdots & 0 & 0 \\
\vdots&\vdots&\vdots&&\vdots&\vdots&&&&\vdots&\vdots&\vdots&& \vdots &\vdots \\
\beta_{s-2}^{1,1}&\beta_{s-3}^{1,1}&\beta_{s-4}^{1,1}&\cdots & \beta_{0}^{1,1}& 0&&\cdots\cdots&&\beta_{s-2}^{3,3}&\beta_{s-3}^{3,3}&\beta_{s-4}^{3,3}&\cdots & \beta_{0}^{3,3} & 0 \\ [5pt]
\beta_{s-1}^{1,1}&\beta_{s-2}^{1,1} &\beta_{s-3}^{1,1} & \cdots& \beta_{1}^{1,1}&\beta_{0}^{1,1} &  &\cdots\cdots &   & \beta_{s-1}^{3,3} &\beta_{s-2}^{3,3} &\beta_{s-3}^{3,3} & \cdots &\beta_{1}^{3,3}& \beta_{0}^{3,3}  \\ [5pt]
\beta_{s}^{1,1}&\beta_{s-1}^{1,1} &\beta_{s-2}^{1,1} & \cdots&\beta_{2}^{1,1}& \beta_{1}^{1,1} &&\cdots\cdots &   & \beta_{s}^{3,3} &\beta_{s-1}^{3,3} & \beta_{s-2}^{3,3}& \cdots & \beta_{2}^{3,3} & \beta_{1}^{3,3}  \\
\vdots&\vdots&\vdots&&\vdots&\vdots&&&&\vdots&\vdots&\vdots&&\vdots&\vdots \\
\beta_{2s+1}^{1,1}&\beta_{2s}^{1,1} &\beta_{2s-1}^{1,1} & \cdots& \beta_{s+3}^{1,1}&\beta_{s+2}^{1,1} &  &\cdots\cdots &   & \beta_{2s+1}^{3,3} &\beta_{2s}^{3,3} &\beta_{2s-1}^{3,3} & \cdots&\beta_{s+3}^{3,3} & \beta_{s+2}^{3,3}  \\ [5pt]
\beta_{2s+2}^{1,1}&\beta_{2s+1}^{1,1} &\beta_{2s}^{1,1} & \cdots&\beta_{s+4}^{1,1} &\beta_{s+3}^{1,1} & &\cdots\cdots &   & \beta_{2s+2}^{3,3} &\beta_{2s+1}^{3,3} &\beta_{2s}^{3,3} & \cdots & \beta_{s+4}^{3,3} &\beta_{s+3}^{3,3}  \\
\vdots&\vdots&\vdots&&\vdots&\vdots&&&&&\vdots&\vdots&&\vdots&\vdots \\
\beta_{4s+1}^{1,1}&\beta_{4s}^{1,1} &\beta_{4s-1}^{1,1} & \cdots& \beta_{3s+3}^{1,1}& \beta_{3s+2}^{1,1} &  &\cdots\cdots &   & \beta_{4s+1}^{3,3} &\beta_{4s}^{3,3} &\beta_{4s-1}^{3,3} & \cdots&\beta_{3s+3}^{3,3} & \beta_{3s+2}^{3,3}  \\ [5pt]
\beta_{4s+2}^{1,1}&\beta_{4s+1}^{1,1} &\beta_{4s}^{1,1} & \cdots&\beta_{3s+4}^{1,1} &\beta_{3s+3}^{1,1} & &\cdots\cdots &  & \beta_{4s+2}^{3,3} &\beta_{4s+1}^{3,3} &\beta_{4s}^{3,3} & \cdots & \beta_{3s+4}^{3,3} & \beta_{3s+3}^{3,3}  \\ [5pt]
0&\beta_{4s+2}^{1,1} &\beta_{4s+1}^{1,1} & \cdots& \beta_{3s+ 5}^{1,1}&\beta_{3s+4}^{1,1}&  &\cdots\cdots &   & 0 &  \beta_{4s+2}^{3,3}&\beta_{4s+1}^{3,3}  & \cdots & \beta_{3s+5}^{3,3} & \beta_{4s+4}^{3,3}  \\ [5pt]
0&0 &\beta_{4s+2}^{1,1} & \cdots& \beta_{3s+6}^{1,1}&\beta_{3s+5}^{1,1}& &\cdots\cdots &   & 0 &  0& \beta_{4s+2}^{3,3}  & \cdots & \beta_{3s+6}^{3,3} & \beta_{3s+5}^{3,3}  \\ [5pt]
0&0 &0 & \cdots& \beta_{3s+7}^{1,1}&\beta_{3s+6}^{1,1}& &\cdots\cdots &   & 0 &  0& 0  & \cdots & \beta_{3s+7}^{3,3} & \beta_{3s+6}^{3,3}  \\
\vdots&\vdots&\vdots&&\vdots&\vdots&&&&\vdots&\vdots&\vdots&&\vdots&\vdots \\
0&0 & 0&\cdots& \beta_{4s+2}^{1,1}& \beta_{4s+1}^{1,1} &   &\cdots\cdots &   & 0   & 0& 0& \cdots& \beta_{4s+2}^{3,3} & \beta_{4s+1}^{3,3}  \\ [5pt]
0&0 & 0&\cdots& 0& \beta_{4s+2}^{1,1} &   &\cdots\cdots &   & 0   &0 & 0& \cdots &0& \beta_{4s+2}^{3,3}  \\
\end{array}\right]$$

\smallskip

Now consider the matrix $A$ in Lemma~\ref{main_lemma}. Specializing the variables $T_l^{(t)}$ to the coefficients $\gamma_{r}^{(j)}$, by the definition of the $\beta$'s one gets that the entries of $A$ specialize to the entries of $\mathcal A$ above.
Since  $A$ has maximal rank then for general choice of $\gamma$'s, so will $\mathcal A$.  But $6s\geq 5s+2 \Leftrightarrow s\geq 2$. Therefore, the linear map defined by $\mathcal A$ is surjective, i.e., the linear system has solution for general values of $\gamma$'s.

\smallskip

The entire argument works in exactly the same way for $f_2,f_3$, hence one has $y^{3s}f_j\in I^{2},$ for $j=1,2,3.$ 

\smallskip

\noindent {\bf Even case:} $d=2s$. The argument is analogous to the one in the odd case, except that one is assuming that $s\geq 5$.
At the end one has to specialize the variables $T_l^{(t)}$ so that the matrix $A$ in Lemma~\ref{main_lemma} (b) specializes to the following content matrix $\mathcal A$:

$$\small\left[\begin{array}{cccccccccccccccccccc}
\beta_{0}^{1,1}& 0& 0& \cdots& 0 & 0&  &\cdots\cdots &  & \beta_{0}^{3,3} & 0& 0& \cdots & 0 & 0  \\ [5pt]
\beta_{1}^{1,1}& \beta_{0}^{1,1}& 0& \cdots& 0 &0&   &\cdots\cdots &   & \beta_{1}^{3,3} & \beta_{0}^{3,3}& 0& \cdots & 0 & 0  \\ [5pt]
\beta_{2}^{1,1}&\beta_{1}^{1,1}&\beta_{0}^{1,1}&\cdots & 0& 0&&\cdots\cdots&&\beta_{2}^{3,3}&\beta_{1}^{3,3}&\beta_{0}^{3,3}&\cdots & 0 & 0 \\
\vdots&\vdots&\vdots&&\vdots&\vdots&&&&\vdots&\vdots&\vdots&& \vdots &\vdots \\
\beta_{s-3}^{1,1}&\beta_{s-4}^{1,1}&\beta_{s-5}^{1,1}&\cdots & \beta_{0}^{1,1}& 0&&\cdots\cdots&&\beta_{s-3}^{3,3}&\beta_{s-4}^{3,3}&\beta_{s-5}^{3,3}&\cdots & \beta_{0}^{3,3} & 0 \\ [5pt]
\beta_{s-2}^{1,1}&\beta_{s-3}^{1,1} &\beta_{s-4}^{1,1} & \cdots& \beta_{1}^{1,1}&\beta_{0}^{1,1} &  &\cdots\cdots &   & \beta_{s-2}^{3,3} &\beta_{s-3}^{3,3} &\beta_{s-4}^{3,3} & \cdots &\beta_{1}^{3,3}& \beta_{0}^{3,3}  \\ [5pt]
\beta_{s-1}^{1,1}&\beta_{s-2}^{1,1} &\beta_{s-3}^{1,1} & \cdots&\beta_{2}^{1,1}& \beta_{1}^{1,1} &&\cdots\cdots &   & \beta_{s-1}^{3,3} &\beta_{s-2}^{3,3} & \beta_{s-3}^{3,3}& \cdots & \beta_{2}^{3,3} & \beta_{1}^{3,3}  \\
\vdots&\vdots&\vdots&&\vdots&\vdots&&&&\vdots&\vdots&\vdots&&\vdots&\vdots \\
\beta_{2s}^{1,1}&\beta_{2s-1}^{1,1} &\beta_{2s-2}^{1,1} & \cdots& \beta_{s+3}^{1,1}&\beta_{s+2}^{1,1} &  &\cdots\cdots &   & \beta_{2s}^{3,3} &\beta_{2s-1}^{3,3} &\beta_{2s-2}^{3,3} & \cdots&\beta_{s+3}^{3,3} & \beta_{s+2}^{3,3}  \\ [5pt]
\beta_{2s+1}^{1,1}&\beta_{2s}^{1,1} &\beta_{2s-1}^{1,1} & \cdots&\beta_{s+4}^{1,1} &\beta_{s+3}^{1,1} & &\cdots\cdots &   & \beta_{2s+1}^{3,3} &\beta_{2s}^{3,3} &\beta_{2s-1}^{3,3} & \cdots & \beta_{s+4}^{3,3} &\beta_{s+3}^{3,3}  \\
\vdots&\vdots&\vdots&&\vdots&\vdots&&&&&\vdots&\vdots&&\vdots&\vdots \\
\beta_{4s-1}^{1,1}&\beta_{4s-2}^{1,1} &\beta_{4s-3}^{1,1} & \cdots& \beta_{3s+2}^{1,1}& \beta_{3s+1}^{1,1} &  &\cdots\cdots &   & \beta_{4s-1}^{3,3} &\beta_{4s-2}^{3,3} &\beta_{4s-3}^{3,3} & \cdots&\beta_{3s+2}^{3,3} & \beta_{3s+1}^{3,3}  \\ [5pt]
\beta_{4s}^{1,1}&\beta_{4s-1}^{1,1} &\beta_{4s-2}^{1,1} & \cdots&\beta_{3s+3}^{1,1} &\beta_{3s+2}^{1,1} & &\cdots\cdots &  & \beta_{4s}^{3,3} &\beta_{4s-1}^{3,3} &\beta_{4s-2}^{3,3} & \cdots & \beta_{3s+3}^{3,3} & \beta_{3s+2}^{3,3}  \\ [5pt]
0&\beta_{4s}^{1,1} &\beta_{4s-1}^{1,1} & \cdots& \beta_{3s+ 4}^{1,1}&\beta_{3s+3}^{1,1}&  &\cdots\cdots &   & 0 &  \beta_{4s}^{3,3}&\beta_{4s-1}^{3,3}  & \cdots & \beta_{3s+4}^{3,3} & \beta_{4s+3}^{3,3}  \\ [5pt]
0&0 &\beta_{4s}^{1,1} & \cdots& \beta_{3s+5}^{1,1}&\beta_{3s+4}^{1,1}& &\cdots\cdots &   & 0 &  0& \beta_{4s}^{3,3}  & \cdots & \beta_{3s+5}^{3,3} & \beta_{3s+4}^{3,3}  \\ [5pt]
0&0 &0 & \cdots& \beta_{3s+6}^{1,1}&\beta_{3s+5}^{1,1}& &\cdots\cdots &   & 0 &  0& 0  & \cdots & \beta_{3s+6}^{3,3} & \beta_{3s+5}^{3,3}  \\
\vdots&\vdots&\vdots&&\vdots&\vdots&&&&\vdots&\vdots&\vdots&&\vdots&\vdots \\
0&0 & 0&\cdots& \beta_{4s}^{1,1}& \beta_{4s-1}^{1,1} &   &\cdots\cdots &   & 0   & 0& 0& \cdots& \beta_{4s}^{3,3} & \beta_{4s-1}^{3,3}  \\ [5pt]
0&0 & 0&\cdots& 0& \beta_{4s}^{1,1} &   &\cdots\cdots &   & 0   &0 & 0& \cdots &0& \beta_{4s}^{3,3} 
\end{array}\right]$$
\qed

\section{The main result}

\begin{Theorem}\label{MAIN}
	Let $I\subset R:=k[x,y]$ denote a codimension $2$ ideal generated by three general forms $f_1,f_2,f_3$ of degree $d\geq 5$.
	Then the Rees algebra of $I$ has depth one.
\end{Theorem}
\demo First assume that, in addition,  $d\geq 10$ when $d$ is even. By Proposition~\ref{annihilators},  the Ratliff--Rush closure of $I$ is strictly larger than $I$. Therefore, by \cite[(1.2)]{HLS} the associated graded ring of $I$ has depth zero and hence the Rees algebra of $I$ has depth $1$.

Thus, it remains to consider the situation where $d$ is even and $5\leq d\leq 9$, i.e., when $d=6$ or $8$.
We will establish below a more general result that may be of interest in itself.
The agenda at this point is to use the Huckaba--Marley test (\cite[Theorem 4.7 (b)]{HM}) for these low degrees.

As a matter of precision, the criterion of Huckaba--Marley is stated for local rings (\cite{HM}). Here, as commonly done, whenever using the local tool, we harmlessly pass to the local ring $k[x,y]_{\fm}$, where $\fm=(x,y)$, and consider the extended ideal $I_{\fm}$.
We will make no distinction in the notation.

Let  $\lambda(\_)$ denote length over $R$, let $J\subset I$ stand for a minimal reduction of $I$.
In the next proposition one shows how to translate the test in our context solely in terms of the degree of the given forms.

\begin{Proposition}\label{aci_general}
Let $I\subset R=k[x,y]$ denote a codimension $2$ ideal
generated by $3$ sufficiently  general forms of degree $d\geq 5$.
Suppose that
\begin{equation}\label{HM-rewritten}
\sum_{\ell=3}^{d-1} \lambda (I^{\ell}/JI^{\ell -1})>
\left\{
\begin{array}{lll}
\frac{d(d-2)}{8}+ \max \left\{\frac{d}{2}-5,0\right\} & \mbox{if $d$ is even}\\[8pt]
\frac{(d+1)(d-1)}{8}-2& \mbox{if $d$ is odd}.
\end{array}
\right.
\end{equation}
Then $\mathcal{R}_R(I)$ has depth one.
\end{Proposition}
\demo
Let $e_1(I)$ stand for the second coefficient
of the Hilbert--Samuel polynomial of $R/I$ in its combinatorial expression.
According to the Huckaba--Marley approach the strict inequality
\begin{equation}\label{HM}
\sum_{\ell\geq 1} \lambda (I^{\ell}/JI^{\ell -1}) > e_1(I),
\end{equation}
implies that the associated graded ring of $I$ has depth zero, hence that $\mathcal{R}_R(I)$ has depth one.

Since the generators of $I$ are general forms of the same degree, one can assume that two of them  will generate a minimal reduction $J\subset I$ -- this is a first weak use of general forms.

We recall the easy general fact that, since $I=(J,f)$, for some $f\in I$, then $I^{\ell}=(JI^{\ell-1}, f^{\ell})$,
for any $\ell\geq 1$.
It follows that
\begin{eqnarray*} I^{\ell}/JI^{\ell -1}&\simeq &(JI^{\ell-1}, f^{\ell})/JI^{\ell -1} \simeq
(f^{\ell})/JI^{\ell -1}\cap (f^{\ell})\\ \nonumber
&\simeq & (f^{\ell})/(JI^{\ell -1}:f^{\ell}) (f^{\ell})\simeq R/JI^{\ell -1}:f^{\ell}, \nonumber
\end{eqnarray*}
hence $\lambda ( I^{\ell}/JI^{\ell -1})=\lambda (R/JI^{\ell-1}:f^{\ell})$.

Now, taking a set of generators of $I^{\ell}$ where $f^{\ell}$ is listed last, then $JI^{\ell-1}:f^{\ell}$ coincides with the ideal $\mathfrak{a}_{\ell}\subset R$ generated by the corresponding syzygy coordinates.
Clearly, $\mathfrak{a}_{\ell}$ is an $(x,y)$-primary ideal.
Thus, $\lambda ( I^{\ell}/JI^{\ell -1})=\lambda (R/\mathfrak{a}_{\ell})$.
It then follows immediately from Lemma~\ref{semi-balanced} that
$$\lambda (I/J)=
\left\{
\begin{array}{lll}
\frac{d^2}{4} & \mbox{if $d$ is even}\\[5pt]
\frac{d^2-1}{4} & \mbox{if $d$ is odd.}
\end{array}
\right.
$$

\medskip

Next, by computing the Hilbert function of $R/I_1(\phi)$ using Lemma~\ref{semi-balanced} again, one has:
$$
\lambda(R/I_1(\phi))=
{{\left \lfloor{\frac{d}{2}}\right \rfloor +1}\choose {2}} +
\left\{
\begin{array}{lll}
\max \{\frac{d}{2}-5,0\} & \mbox{if $d$ is even} \\[5pt]
\;\;\left \lfloor{\frac{d}{2}}\right \rfloor -2 & \mbox{if $d$ is odd.}
\end{array}
\right.
$$
(A second weak use of general forms.)

On the other hand, by \cite[Proposition 3.12]{syl3} one has $\lambda(I^2/JI)=\lambda(I/J)-\lambda(R/I_1(\phi))$.
Therefore, it obtains
$$
\lambda(I^2/JI)=
{{\left \lceil{\frac{d}{2}}\right \rceil}\choose {2}} -
\left\{
\begin{array}{lll}
\max \{\frac{d}{2}-5,0\} & \mbox{if $d$ is even} \\[5pt]
\;\;\left(\left \lfloor{\frac{d}{2}}\right \rfloor -2\right) & \mbox{if $d$ is odd.}
\end{array}
\right.
$$

\smallskip

We next claim that the rational map $\pp^1\dasharrow \pp^2$ defined by the generators of $I$ is birational.
Indeed, one can reduce the problem to the affine situation by setting $t:=x/y$ in the usual way after dividing all
three generators of $I$ by $y^{d}$  and then taking the fractions with same denominator (one of the three).
This way we get a rational map $\mathbb{A}^1\dasharrow \mathbb{A}^2$ defined by rational functions
$F_1(t)/F_3(t), F_2(t)/F_3(t)$, where $\deg F_i(t)=d$.
Since the involved terms are general polynomials in $t$ ($k$ being infinite), for a general
point $(a_1,a_2)\in \mathbb{A}^2$ the system of equations
$F_1(t)/F_3(t)=a_1, F_2(t)/F_3(t)=a_2$
admits exactly one solution, including multiplicity (algebraically: $t$ is a rational fraction
in $F_1(t)/F_3(t), F_2(t)/F_3(t)$).

(Yet another use of general forms.)

\smallskip

As a consequence, one has $e_1(I)=\frac{1}{2} (d^2-d)= {{d}\choose {2}}$ and $\deg(R/I)=d$
(see \cite[Proposition 3.3]{syl2}).
It also follows that the reduction number of $I$ is $ d-1$.

Thus, the Huckaba--Marley criterion requires that for $d\geq 5$ one have
$$\sum_{\ell=1}^{d-1} \lambda_{\ell}(I^{\ell}/JI^{\ell-1}) > {{d}\choose {2}}.$$
Subtracting the sum $\lambda (I/J)+\lambda (I^2/JI)$ as computed above, it suffices to guarantee the inequalities 
$$\sum_{\ell=3}^{d-1} \lambda_{\ell}(I^{\ell}/JI^{\ell-1})>
\left\{
\begin{array}{lll}
\frac{d(d-2)}{8}+ \max \left\{\frac{d}{2}-5,0\right\} & \mbox{if $d$ is even}\\[8pt]
\frac{(d+1)(d-1)}{8}-2& \mbox{if $d$ is odd},
\end{array}
\right.
$$
as was to be shown.
\qed

\medskip

We now apply the above to conclude the proof of Theorem~\ref{MAIN}.

For simplicity, fixing the minimal reduction $J\subset I$, we set $\lambda_{\ell}:=\lambda (I^{\ell}/JI^{\ell -1})$.
As one easily realizes, for $5\leq d\leq 11$ it will suffice to compute $\lambda_3$ and bound below each of the remaining $\lambda_{\ell}$'s ($\ell\geq 4$) by $1$.
(We remark that for higher values of $d$ (e.g., $d=13$) more precise bounds for $\lambda_4$, etc., may be required for the conclusion.) 

As noted in the proof of Proposition~\ref{HM-rewritten}, finding $\lambda_3$ depends on determining a  minimal set of generators of the ideal $\mathfrak{a}_3$ generated by the last coordinates of the syzygies of $I^3$.
Computing with \cite{Macaulay1} by employing random forms and using the established values of $\lambda_1, \lambda_2$ in the proof of the proposition, one finds: 

\smallskip



$d=6$: $\mathfrak{a}_3 = (x,y)^2 \Rightarrow \lambda_3+\lambda_4+\lambda_5\geq 5$; since $\lambda_1=d^2/4=9$ and $\lambda_2\geq \lambda_3$ the total sum is at least $9+3+5=17>15={{6}\choose {2}}$




$d=8$: $\mathfrak{a}_3 = (x,y)^2 \Rightarrow \lambda_3=3 \Rightarrow \lambda_3+\cdots +\lambda_7\geq 7$; since $\lambda_1=64/4=16$ and $\lambda_2=6$ for $d=8$, we find the total sum is at least $29>28={{8}\choose {2}}$







\medskip

This completes the proof of Theorem~\ref{MAIN}.
\qed

\begin{Remark}\rm 
A purely theoretical argument, without the computer, even for $d=6$, is quite subtle.
Here for the failure of the almost Cohen--Macaulay property we need to show the bound $\lambda_3+\lambda_4+\lambda_5 > 3$, so it suffices to show that
$\lambda_3\geq 2$ since the reduction number is $5$.
Equivalently, we are to show that $I^3$ has at most one linear syzygy because then, since $\mathfrak{a}_i$ is always an ideal of finite colength, $\lambda_3\geq 2$.
Contradicting this takes us to a long discussion about the degrees of the syzygies of $I^3$ which eventually abuts 
at the following:
since $I^3$ is a perfect ideal with $10$ generators and generated in degree $18$, its minimal presentation matrix has
$7$ columns whose degrees $r_3\leq\cdots \leq r_9$ add up to $16$. By elementary column operations with pivot
the last coordinates of the two assumed linear syzygies, we may assume that the last entry of any of the
other $7$ columns are zero.
Then they are syzygies of $JI^2$.
Therefore, $r_9\geq \cdots \geq r_3\geq 2$. On the other hand, no minimal syzygy has degree $\geq 4$
since the presentation matrix of $I$ is in degree $3$. To add up to $16$, the only way is
$r_3=\cdots =r_7=2$, $r_8=r_9=3$.
This structure of the Hilbert--Burch matrix of $I^3$ may be attained
if the entries of the Hilbert--Burch matrix of $I$ are not general - see Example~\ref{aci_nongeneral} (b) below.

The theoretical side of the discussion for $d\geq 7$ eludes the eye.
\end{Remark}

The following examples show that, already for $d=6$, the statement of Theorem~\ref{MAIN} or the hypothesis of Theorem~\ref{aci_general} are																																																																																																																		 no longer
true  if the generators of the ideal are not general enough,
even when the degrees of the syzygies are as in Lemma~\ref{semi-balanced}.

\begin{Example}\label{aci_nongeneral}\rm
(a) The first example keeps some of the properties above:
1) The associated rational map is birational (onto the image), hence $e_1(I)={{6}\choose {2}}=15$
and the reduction number red$_J(I)= {\rm edeg}(I)-1=2.3-1=5$;
(2) $\lambda(I/J)=9$, but $\lambda(I^2/JI)=2$ instead, a degenerate value;
(3) The claim above that $\lambda_{\ell}={{s_{\ell}+1}\choose {2}}$ fails for $\ell=3$ as here $\lambda_3=2$.

The Hilbert--Burch matrix of this example is
$$\left(
\begin{array}{cc}
-xy^2   &   -y^3\\
x^3  &   0\\
y^3   &   x^3
\end{array}
\right)
$$
A minimal reduction is $(x^6, x^4y^2-y^6)$.
Note that a linear syzygy appears all too soon among the forms $x^6 (x^4y^2-y^6),\, x^6. x^3y^3,\,(x^3y^3)^2$.
One gets here
$$\sum_{\ell=1}^5 \lambda_{\ell}=9 + 2 + 2 +1 +1 = 15.$$

(b) The second example has Hilbert--Burch matrix
$$\left(
\begin{array}{cc}
x^3   &   xy^2\\
x^2y  &   x^3-y^3\\
y^3   &   -x^2y+xy^2
\end{array}
\right)
$$
\end{Example}
The subideal $J:=(x^6-2x^3y^3, x^4y^2-y^6)$ is a minimal reduction.
This example is really on the edge since we can check that:
(1) The associated rational map is again birational, hence $e_1(I)=15$
and  red$_J(I)=5$;
(2) The lengths $\lambda(I/J), \lambda(I^2/JI)$ are as stated in the above preliminaries;
(3) The expected value $\lambda_{\ell}={{s_{\ell}}\choose {2}}$ for $\ell\geq 3$ holds here with $s_3=1$ (i.e., $JI^2:I^3=(x,y)$).

Since $\lambda_3=1$ then $\lambda_5=\lambda_4=1$, where red$_J(I)=5$,
thus yielding
$$\sum_{\ell=1}^5 \lambda_{\ell}=9 + 3 + 1 +1 +1 = 15.$$
Therefore, in both examples the Rees algebra of $I$ has depth $\geq 2$, i.e., it is almost Cohen--Macaulay.



\end{document}